\documentclass[12pt]{article}

\usepackage{amscd,amsmath, amssymb, fancyhdr, mathbbol}

\usepackage[backref=page]{hyperref}
\renewcommand*{\backrefalt}[4]{%
	\ifcase #1 (Not cited.)%
	\or        (Cited on page~#2.)%
	\else      (Cited on pages~#2.)%
	\fi}

\hypersetup{
	colorlinks   = true,
	citecolor    = magenta,
	linkcolor    = blue,
	urlcolor     = magenta	
}

\numberwithin{equation}{section}

\newcommand{\version}{version 1.1,\ \ Sep. 10, 2026}

\def\eqref#1{(\ref{#1})}

\newcommand{\arrow}{{\:\longrightarrow\:}}
\newcommand{\Z}{{\Bbb Z}}
\def\C{{\Bbb C}}

\newcommand{\Q}{{\Bbb Q}}

\newcommand{\6}{\partial}
\def\1{\sqrt{-1}\:}
\newcommand{\restrict}[1]{{\left|_{{\phantom{|}\!\!}_{#1}}\right.}}
\newcommand{\cntrct}                
{\hspace{2pt}\raisebox{1pt}{\text{$\lrcorner$}}\hspace{2pt}}

\newcommand{\calo}{{\cal O}}

\renewcommand{\tilde}{\widetilde}
\renewcommand{\bar}{\overline}
\renewcommand{\phi}{\varphi}
\renewcommand{\epsilon}{\varepsilon}
\renewcommand{\geq}{\geqslant}

\newcommand{\Tot}{\operatorname{Tot}}

\newcommand{\Pic}{\operatorname{Pic}}
\newcommand{\Ext}{\operatorname{Ext}}
\newcommand{\Hom}{\operatorname{Hom}}

\newcommand{\const}{\operatorname{\sf const}}
\newcommand{\slope}{\operatorname{slope}}

\newcounter{Mycounter}[section]
\newcounter{lemma}[section]
\renewcommand{\thelemma}{{Lemma \thesection.\arabic{lemma}}}
\newcommand{\lemma}{%
    \setcounter{lemma}{\value{Mycounter}}
    \refstepcounter{lemma}
    \stepcounter{Mycounter}
    {\noindent \bf \thelemma:\ }}

\newcounter{claim}[section]
\renewcommand{\theclaim}{{Claim \thesection.\arabic{claim}}}
\newcommand{\claim}{%
    \setcounter{claim}{\value{Mycounter}}
    \refstepcounter{claim}
    \stepcounter{Mycounter}
    {\noindent \bf \theclaim:\ }}

\newcounter{sublemma}[section]
\newcounter{corollary}[section]
\renewcommand{\thecorollary}{{Corollary \thesection.\arabic{corollary}}}
\newcommand{\corollary}{%
    \setcounter{corollary}{\value{Mycounter}}
    \refstepcounter{corollary}
    \stepcounter{Mycounter}
    {\noindent \bf \thecorollary:\ }}

\newcounter{theorem}[section]
\renewcommand{\thetheorem}{{Theorem \thesection.\arabic{theorem}}}
\newcommand{\theorem}{%
    \setcounter{theorem}{\value{Mycounter}}
    \refstepcounter{theorem}
    \stepcounter{Mycounter}
    {\noindent \bf \thetheorem:\ }}

\newcounter{conjecture}[section]
\renewcommand{\theconjecture}{{Conjecture \thesection.\arabic{conjecture}}}
\newcommand{\conjecture}{%
    \setcounter{conjecture}{\value{Mycounter}}
    \refstepcounter{conjecture}
    \stepcounter{Mycounter}
    {\noindent \bf \theconjecture:\ }}

\newcounter{proposition}[section]
\renewcommand{\theproposition}
      {{Proposition \thesection.\arabic{proposition}}}
\newcommand{\proposition}{%
    \setcounter{proposition}{\value{Mycounter}}
    \refstepcounter{proposition}
    \stepcounter{Mycounter}
    {\noindent \bf \theproposition:\ }}

\newcounter{definition}[section]
\renewcommand{\thedefinition}
      {{Definition~\thesection.\arabic{definition}}}
\newcommand{\definition}{%
    \setcounter{definition}{\value{Mycounter}}
    \refstepcounter{definition}
    \stepcounter{Mycounter}
    {\noindent \bf \thedefinition:\ }}

\newcounter{example}[section]
\newcounter{remark}[section]
\renewcommand{\theremark}{{Remark \thesection.\arabic{remark}}}
\newcommand{\remark}{%
    \setcounter{remark}{\value{Mycounter}}
    \refstepcounter{remark}
    \stepcounter{Mycounter}
    {\noindent \bf \theremark:\ }}

\newcounter{problem}[section]
\newcounter{question}[section]
\newcommand{\proof}{\noindent{\bf Proof:\ }}

\newcommand{\pstep}{{\bf Proof. Step 1: \ }}

\makeatletter

\@addtoreset{equation}{section} \@addtoreset{footnote}{section}
\makeatother

\def\blacksquare{\hbox{\vrule width 5pt height 5pt depth 0pt}}
\def\endproof{\blacksquare}

\begin{document}
\begin{center}
{\LARGE\bf
Lagrangian fibrations on hyperk\"ahler manifolds have no multiple
fibers in codimension one
\\[4mm]
}

Ljudmila Kamenova\footnote{Partially supported 
by award SFI-MPS-TSM-00013537 from the Simons Foundation International}, 
Misha Verbitsky\footnote{Partially supported 
by FAPERJ SEI-260003/000410/2023 and CNPq - Process 310952/2021-2. 

{\bf 2010 Mathematics Subject
Classification: 53C26, 14J42} }

\end{center}
{\small \hspace{0.10\linewidth}
\begin{minipage}[t]{0.85\linewidth}
{\bf Abstract.} 
Let $f:\; M \to {\mathbb C}P^n$ be 
a Lagrangian fibration on a compact hyperkahler
manifold. We prove that $f$ has no multiple fibers in
codimension 1. We prove that this condition is
equivalent to the primitivity of the fundamental
class of $\pi^{-1}(H)$, where $H \subset {\mathbb C}P^n$
denotes a hyperplane divisor. 
\end{minipage}
}

\tableofcontents


\section{Introduction} 


\subsection{ Multiplicity of a fiber in codimension 1}

In \cite[Theorem 1.1]{_Hwang_Oguiso:multiple_}, 
J.-M. Hwang and K. Oguiso
classified the general singular fibers of 
proper, holomorphically Lagrangian fibrations.
Their analysis did not use the hyperk\"ahler
structure on its total space; indeed, they
were interested in Lagrangian fibrations
over an open ball, with fibers of Fujiki 
class C. This classification is most remarkable
because all the cases listed in 
\cite[Theorem 1.1]{_Hwang_Oguiso:multiple_}
were realized by an explicit construction.
Interestingly enough, some of the examples
constructed by Hwang and Oguiso have multiple fibers
in codimension 1 
(we introduce the notion of multiplicity of 
fibers in \ref{_mult_Definition_} below).

In contrast to the picture presented by Hwang and Oguiso,
the set of possible codimension 1 singular fibers for 
Lagrangian fibrations on compact hyperk\"ahler manifolds
is not known: not all examples in \cite[Theorem 1.1]{_Hwang_Oguiso:multiple_}
can be realized on compact manifolds.

In the present paper, we show that the multiple
fibers cannot be realized, that is, any Lagrangian
fibration from a compact hyperk\"ahler manifold to 
a smooth base is non-multiple in codimension 1.

The writing of this paper took us almost 2 years. Originally,
we wanted to sharpen the results of \cite{_KV:primitive_}
by expressing the multiplicity of the codimension 1 fibers
in terms of the local monodromy of the 
primitive root of $\pi^*(\calo(1))$ 
(\ref{_flat_conn_on_restriction_Definition_}).
By the second year, we realized that these arguments
can be used to show that the local monodromy is trivial,
and there are no multiple fibers of codimension 1.

A few months later, Y. Kim, K. Oguiso, E. Shinder,  and
discovered another proof of this result, based on Koll\'ar
vanishing theorems (\cite{_KOS:obstructions_}). 
They were very kind to suggest us
to publish our result simultaneously, but we were late
with our text for personal reasons. At the end,
their paper appeared on August 10, 2026. However, 
the slides from a talk about an early (and very 
imprecise) version of this proof have been circulated
since December 2025 (\cite{_Verbitsky:talk_Mult_}).

\subsection{ Main results}
\label{_main_Subsection_}

\definition\label{_mult_Definition_}
Let  $\pi:\; M \arrow X$ be a proper holomorphic
map, $x\in X$ a point, $F_x:= \pi^{-1}(x)$ 
and $F_i$ its irreducible component, with 
(scheme-theoretical) multiplicity $\mu_i$.
Denote the greatest common divisor of $\mu_i$
by $\mu$. The number $\mu$ is called {\bf the
multiplicity} of the fiber. A fiber is {\bf multiple} if $\mu >1$.
A fiber $F_x$ is {\bf reduced} 
if $\mu_i=1$ for all $i$.
A fiber $F_x$  {\bf has a  reduced component} 
if $\mu_i=1$ for at least one $i$.

\hfill

The following theorem was proven in
\cite{_CKV:min_mult_}.

\hfill

\theorem 
Let $\pi:\; M \arrow X$ be an abelian fibration
(that is, a fibration with general fiber a complex torus),
and assume that $M$ is K\"ahler. Let $\mu_i$
be the multiplicity of irreducible components
of $\pi^{-1}(z)$ for some $z\in X$,
and $gcd(\mu_i)$ their greatest common divisor. Then $gcd(\mu_i)=\min \mu_i$.

\endproof

\hfill

\definition
Let $\pi:\; M \arrow X$ be an abelian fibration of
complex manifolds
and $D\subset X$ its set of critical values, which
is known as {\bf  the discriminant}, or
{\bf  the discriminant divisor}.
We say that $\pi$ {\bf has no multiple fibers
in codimension 1} if  the fiber $\pi^{-1}(x)$ 
has a component with multiplicity 1 for all
$x$ outside of a codimension 2 subvariety.

\hfill

The following theorem is well-known (\cite{_BHPV_}, 
\cite[Chapter XI, Proposition 1.6]{_Huybrechts:K3_}).

\hfill

\theorem
Let $\pi:\; M \to X$ be an elliptic fibration on a K3 surface.
Then $\pi$ has no multiple fibers.
\endproof

\hfill

The main result of this paper
is a generalization of this theorem.

\hfill

\theorem\label{_main_intro_Theorem_}
Let $\pi:\; M \to \C P^n$ be a
Lagrangian fibration on a hyperk\"ahler manifold.
Then $\pi$ has no multiple fibers
in codimension 1.

\proof \ref{_no_mult_main_end_Theorem_}. \endproof

\hfill

We also prove the following theorem, which is a 
stronger form of the one proved in \cite{_KV:primitive_}.

\hfill

\theorem
Let $\pi:\; M \to \C P^n$ be a
Lagrangian fibration on a hyperk\"ahler manifold,
and $\kappa\in H^2(M, \Z)$ the fundamental class
of $\pi^{-1}(H)$, where $H\subset \C P^n$ is the hyperplane
divisor. Then $\kappa$ is primitive, that is,
not divisible by an integer $n>1$ in $H^2(M, \Z)$.

\proof \ref{_primitivity_Theorem_}. \endproof

\hfill

\remark
In \cite{_KV:primitive_}, this result was proven
assuming $\pi$ has no multiple fibers.
However, a Lagrangian fibration
might have multiple fibers in codimension 2, as
I. Hellmann (\cite[Section 2.2]{_Hellmann:cone_}) has shown.

\hfill

\remark
Using \cite[Theorem 2.10]{_KL:epiga_non-hyperbolic_},
most of the arguments we use can be extended to the
case where $M$ is a  primitive symplectic variety,
and $\pi:\; M \to X$ a Lagrangian fibration 
with normal base admitting a K\"ahler form on
its smooth locus. We expect \ref{_main_intro_Theorem_}
to be true in this setting as well.

\hfill

\remark
In \cite[Theorem 4.12]{_K_Lu:dominability_}, L. Kamenova and S. Lu
prove that for any hyperk\"ahler manifold $M$ admitting a
Lagrangian fibration with no multiple fibers in codimension 1
there exists a holomorphic map $\C^{2n}\to M$
which is locally biholomorphic at some point.
In the present paper we show that this assumption
holds automatically.


\section{Hyperk\"ahler manifolds} 


We start by reminding the reader the basic definitions
and results from the theory of hyperk\"ahler manifolds
(\cite{_Beauville_,_Fujiki:HK_,_Besse:Einst_Manifo_,_Bogomolov:defo_}).
This section contains the preliminary material, 
the reader should feel free to skip it.

\subsection{Hyperk\"ahler structures and Calabi-Yau theorem}

\definition
A {\bf hyperk\"ahler structure} on a manifold $M$
is a Riemannian structure $g$ and a triple of complex
structures $I,J,K$, satisfying quaternionic relations
$I\circ J = - J \circ I =K$, such that $g$ is K\"ahler
for $I,J,K$.

\hfill

\remark A hyperk\"ahler manifold  has three symplectic forms\\
$\omega_I:=  g(I\cdot, \cdot)$, $\omega_J:=  g(J\cdot, \cdot)$,
$\omega_K:=  g(K\cdot, \cdot)$.

\hfill

\remark
This is equivalent to $\nabla I=\nabla J = \nabla K=0$:
the parallel translation along the connection preserves $I, J,K$.

\hfill

\definition  A holomorphically symplectic manifold 
is a complex manifold equipped with non-degenerate, holomorphic
$(2,0)$-form.

\hfill

\remark 
Hyperk\"ahler manifolds are holomorphically symplectic.
Indeed, $\Omega:=\omega_J+\1\omega_K$ is a holomorphic symplectic
form on $(M,I)$.

\hfill

\theorem (Calabi-Yau) 
A compact, K\"ahler, holomorphically symplectic manifold
 admits a unique hyperk\"ahler metric in any K\"ahler class.

\subsection{Hyperk\"ahler manifolds of maximal holonomy}

\definition A hyperk\"ahler manifold $M$ is called
{\bf of maximal holonomy}, or {\bf IHS}
if $\pi_1(M)=0$, $H^{2,0}(M)=\C$. 

\hfill

As follows from Bogomolov's decomposition theorem
(\cite{_Bogomolov:decompo_}), any 
hyperk\"ahler manifold admits a finite covering
which is a product of a torus and several 
hyperk\"ahler manifolds of maximal holonomy.
Maximal holonomy implies the following 
remarkable topological result, due to A. Fujiki.

\hfill

\theorem (A. Fujiki, \cite{_Fujiki:HK_})\\
Let $\eta\in H^2(M)$, and $\dim M=2n$, where $M$ is
hyperk\"ahler. Then 
\begin{equation}\label{_Fujiki_Equation_}
\int_M \eta^{2n}=cq(\eta,\eta)^n,
\end{equation}
for some primitive integer quadratic form $q$ on $H^2(M,\Z)$,
and $c>0$ a rational number.

\hfill

\definition
This form is called
{\bf Bogomolov-Beauville-Fujiki form}. It is defined
by the Fujiki's relation \eqref{_Fujiki_Equation_}
uniquely, up to a sign. The sign is determined
from the following formula (Bogomolov, Beauville)
\begin{align*}  \lambda q(\eta,\eta) &=
   \int_X \eta\wedge\eta  \wedge \Omega^{n-1}
   \wedge \bar \Omega^{n-1} -\\
 &-\frac {n-1}{2n}\left(\int_X \eta \wedge \Omega^{n-1}\wedge \bar
   \Omega^{n}\right) 
\left(\int_X \eta \wedge \Omega^{n}\wedge \bar \Omega^{n-1}\right)
\end{align*}
where $\Omega$ is the holomorphic symplectic form, and 
$\lambda>0$.

\hfill

  For the rest of this paper,
a hyperk\"ahler manifold
is a compact, K\"ahler, holomorphically symplectic manifold
of maximal holonomy.

\subsection{Lagrangian fibrations}

In this subsection,
we define a Lagrangian fibration on a hyperk\"ahler
manifold and cite basic (and foundational) results 
of D. Matsushita and J.-M. Hwang.

\hfill

\definition
Let $M$ be a holomorphically symplectic manifold
and  $\pi:\; M \to B$ a holomorphic map. We say that
$\pi$ is {\bf a Lagrangian fibration} if any smooth fiber
$F$ of $\pi$ is a connected holomorphic Lagrangian submanifold of $M$ 
(that is, its dimension is $\frac 1 2 \dim M$ and the holomorphically
symplectic form vanishes on $F$). 
Its {\bf discriminant} is the set of its critical values;
by Sard lemma and Remmert theorem, the discriminant is
a proper complex subvariety in $B$. 
J.-M. Hwang and K. Oguiso prove that the
discriminant of a Lagrangian fibration
is always a divisor 
(\cite[Proposition 3.1]{_Hwang_Oguiso:characteristic_}). 

\hfill

\theorem (D. Matsushita, \cite{_Matsushita:fibred_})\\
Let $M$ be a hyperk\"ahler manifold of maximal holonomy,
and $\pi:\; M \arrow B$ a surjective holomorphic map,
with $0 < \dim B < \dim M$. Then $\pi$ is
a Lagrangian fibration.

\hfill

\theorem (J.-M. Hwang, \cite{_Hwang:CP^n_})\\
In these assumptions,  $B$ is biholomorphic to $\C P^n$ 
when it is smooth.

\hfill

\conjecture 
The base $B$ 
of a Lagrangian fibration
is biholomorphic to $\C P^n$ when it is normal.

\hfill

\remark In \cite{_CMS-B:AdvStud_}, Cho, Miyaoka, Shepherd-Barron 
stated this conjecture as a theorem, but the proof
had gaps (\cite{_Kebekus_}). When $\dim_\C M=4$, it was proven
in \cite{_Hu_Chenyang:4folds_}.

\hfill

\theorem (D. Matsushita, \cite{_Matsushita:CP^n_})\\
Let $M$ be a hyperk\"ahler manifold of maximal holonomy,
and $\pi:\; M \arrow X$ a Lagrangian fibration, with $X$ normal.
Then $H^*(X, \Q)\cong H^*(\C P^n, \Q)$.

\hfill

\remark\label{_AL_Remark_}
General fibers of $\pi$ are abelian varieties (projective
complex tori), by Arnold-Liouville theorem 
\cite{_Evans:Lagrangian_}. Conversely, as shown by Hwang-Weiss, 
any Lagrangian complex torus in $M$ is a fiber of a Lagrangian fibration
(\cite{_Hwang_Weiss_}).


\section{The ETMDPS vanishing theorem and its applications}


\subsection{Line bundles and the Riemann-Roch formula}

Recall that a class $\eta\in H_k(M,\Z)$ is called
{\bf primitive} if it is not divisible, that
is, there is no $\eta'\in H_k(M,\Z)$ such that
$\eta= r\eta'$, with $r\in \Z$, $|r|\geq 2$.
A line bundle $L$ is {\bf primitive} if $c_1(L)\in H^2(M, \Z)$
is a primitive class.

\hfill

Let $\pi:\; M \to \C P^n$ be a Lagrangian fibration,
and $\calo(1)$ the hyperplane line bundle on $\C P^n$.
In \cite{_KV:primitive_}, we show that when all fibers of $\pi$
are reduced, the pullback $\pi^*\calo(1)$ is primitive.
In Section \ref{_primitivity_Section_}, we extend this argument to fibrations
without multiple fibers in codimension 1.
The argument in \cite{_KV:primitive_}
was based on the following version of Riemann-Roch 
formula, due to D. Huybrechts.

\hfill

\theorem \label{_Huybrechs_chi_Theorem_}
Let $M$ be a hyperk\"ahler manifold, $\dim_C M=2n$
and $L$ a holomorphic line bundle.
Then $\chi (L) = \sum a_i q(c_1(L))^i$, where the coefficients $a_i$ are 
rational constants depending on the topology of $M$. 

\proof
\cite[Section 1.11]{_Huybrechts:basic_}.
\endproof

\hfill

We will only use the following corollary of this theorem.

\hfill

\corollary\label{_chi(L)_Corollary_}
Let $\pi:\; M \to X$ be a Lagrangian 
fibration on a hyperk\"ahler manifold, 
and $L$ a line bundle on $M$ such that
$L^{\otimes d}= \pi^*(L_0)$, for some $d \in \Z$
where $L_0$ is a line bundle on
$X$. Then $\chi(L)=n+1$, where $2n= \dim_\C M$.

\hfill

\proof
Since $\chi(\calo_M)=n+1$, 
\ref{_Huybrechs_chi_Theorem_} implies that
$\chi(L)=n+1$ for any $L$ such that
$q(c_1(L), c_1(L))=0$. On the other hand,
by Fujiki formula, $q(\eta,\eta)^n= \const \int_M \eta^{2n}$,
and this integral manifestly vanishes for any 
2-form $\eta$ obtained as a pullback from $X$.
\endproof

\hfill

The following topological observation can be used to
avoid the ambiguity when taking $L$ such that
$L^{\otimes d}=\pi^*(\calo_{\C P^n}(1))$.

\hfill

\claim\label{_torsion_free_Claim_}
Let $M$ be a hyperk\"ahler manifold of maximal holonomy.
Then $H^2(M)$ is torsion-free.

\hfill

\proof The universal coefficients formula gives
the exact sequence: 
\[
  0 \to \Ext_\Z^1(H_1(X; \Z), \Z) \to H^2(X; \Z) \to
    \Hom_\Z(H_2(X; \Z), \Z)\to 0.
\]
Since $H_1(X, \Z)=0$ for a maximal holonomy hyperk\"ahler manifold, 
this gives an isomorphism $H^2(X;\Z)=\Hom_\Z(H_2(X; \Z), \Z)$,
hence the torsion vanishes.
\endproof

\subsection{The ETMDPS vanishing theorem}

\definition
A real $(1,1)$-form $\eta$ on a complex manifold $M$
is called {\bf semipositive} if $\eta(x, Ix) \geq 0$
for all real tangent vectors $x$.

\hfill

The following theorem was rediscovered several times during
1990-ies (\cite{_Enoki:semipositive_,_Mourougane_,_Takegoshi_}). 
Its most general form (which we
do not use) is due to Demailly, Peternell and Schneider
\cite[Corollary 2.1.2]{_DPS_}. We call this theorem
``the ETMDPS vanishing theorem'', after
Enoki, Takegoshi, Morougane, Demailly, Peternell and Schneider.

\hfill

\theorem\label{_DPS_Theorem_}
Let $(M, I, \omega)$ be a compact K\"ahler manifold, $\dim_\C M=n$, 
$K$ its canonical bundle,
and $L$ a holomorphic line bundle on $M$ equipped with a 
Hermitian metric $h$. Assume that the 
curvature $\Theta$ of $L$ is a semipositive form on $M$.
Then 
the wedge multiplication operator 
$\eta \arrow \omega^i \wedge \eta$
induces a surjective map
\[
H^0(\Omega^{n-i}M\otimes  L)\stackrel {\omega^i \wedge
\cdot}\arrow H^i(K  \otimes  L).
\]
Here $\omega$ is considered as an element in $H^1(\Omega^1 M)$,
and multiplication by $\omega$ maps $H^k(\Omega^{n-l}M\otimes L)$
to $H^{k+1}(\Omega^{n-l+1}M\otimes L)$.
\endproof

\hfill

\corollary\label{_n+1+H_0_Omega^iotimesL_Corollary_}
Let $L$ be a line bundle on a hyperk\"ahler manifold $M$,
admitting a connection with semi-positive curvature
and satisfying $q(c_1(L),c_1(L))=0$. Then 
$\sum_i \dim H^0(M,\Omega^iM\otimes L)\geq n+1$.

\hfill

\proof By 
the ETMDPS vanishing theorem (\ref{_DPS_Theorem_}),
\[ \sum_i \dim H^0(M,\Omega^iM\otimes L)\geq \chi(L), \]
and \ref{_chi(L)_Corollary_} implies that $\chi(L) =n+1$.
\endproof

\hfill

To apply this result, we 
state several elementary algebro-geometric observations.

\hfill

\lemma\label{_vanishing_filtration_Lemma_}
Let $B$ be a vector bundle equipped with a filtration
$0=B_0\subset B_1 \subset ... \subset B_k=B$.
Assume that $H^0(B_i/B_{i-1})=0$. Then $H^0(B)=0$.

\hfill

\proof The proof is based on induction on $k$ 
and uses the long exact sequence of cohomology
associated with the short exact sequences
$0\to B_{k-1} \to B_k \to B_k/B_{k-1}\to 0$.
\endproof

\hfill

\theorem\label{_trivial_on_fibers_Theorem_}
Let $M$ be a hyperk\"ahler manifold admitting a Lagrangian fibration 
$\pi:\; M \rightarrow X$, and $H$ a line bundle on $X$.
Let $L$ be a line bundle such that $L^{\otimes k}=\pi^* H$.
Then $L$ is trivial on all smooth fibers of $\pi$. 

\hfill

\pstep 
Let $F$ be a smooth fiber of $\pi$, which is an abelian variety
by Arnol'd-Liouville (\ref{_AL_Remark_}). 
Then $T^*M\restrict F$ is an extension
of a trivial bundle $TF$ with another trivial bundle $NF= T^* F$.
For any non-trivial line bundle $L\in \Pic_0(F)$, we have
$H^0(L \otimes TF)=0$ and $H^0(L \otimes NF)=0$,
which implies that  $H^0(L \otimes  T^*M\restrict F)=0$.
Similarly, one obtains $H^0(L \otimes \Lambda^k T^*M\restrict F)=0$
(\ref{_vanishing_filtration_Lemma_}).

\hfill

{\bf Step 2:}
Unless $L$ is trivial on $F$, we 
have $H^0(L \otimes \Lambda^k T^*M\restrict F)=0$,
which implies $H^0(L \otimes \Lambda^* M)=0$.
By the ETMDPS vanishing theorem
this gives $H^i(L)=0$, hence $\chi(L)=0$,
contradicting the formula $\chi(L)=n+1$ (\ref{_chi(L)_Corollary_}).
\endproof

\section{Primitivity of line bundles}
\label{_primitivity_Section_}

\subsection{Flat connection on the primitive root of $\pi^*(\calo(1))$}

Let $\pi:\; M \to \C P^n$ be a Lagrangian fibration
on a hyperk\"ahler manifold.
In \cite{_KV:primitive_}, we have shown that
the bundle $\pi^*\calo(1)$  is primitive if all
fibers of $\pi$ are reduced.
In this subsection, we give a new version of this
argument, proving the following theorem.

\hfill

\theorem\label{_primitivity_Theorem_}
Let $\pi:\; M \to \C P^n$ be a Lagrangian fibration
on a hyperk\"ahler manifold, without  multiple
fibers in codimension 1. Then $\pi^*\calo(1)$ 
is a primitive line bundle on $M$. 

\proof See the end of this subsection. \endproof

\hfill

Let $\pi:\; M \to \C P^n$ be a Lagrangian fibration,
$H \subset \C P^n$ a hyperplane section, and
$M_0:= M \backslash \pi^{-1}(H)=\pi^{-1}(\C^n)$. Consider a 
line bundle $L$ on $M$ such that $L^{\otimes d}=\pi^*(\calo_{\C P^n}(1))$.
We study the natural flat connection
on $L\restrict {M_0}$ and compute its monodromy
to prove that $\pi^*\calo(1)$ is primitive.

Recall that any Hermitian holomorphic bundle $(B, h)$
admits a unique connection $\nabla$ such that
$\nabla^{0,1}=\bar\6$ is its holomorphic structure operator,
and $\nabla(h)=0$. This connection is called {\bf the Chern 
connection}.

\hfill

\proposition\label{_flat_connection_Proposition_}
Let $M$ be a hyperk\"ahler manifold admitting a Lagrangian fibration 
$\pi:\; M \rightarrow X$, and $L_0$ a line bundle on $X$.
Let $L$ be a line bundle such that $L^{\otimes k}=\pi^* L_0$.
Then $L$ admits a Chern connection $\nabla$ which is flat
on each restriction $L\restrict F$ to the fiber $F$ of $\pi$.
Moreover, for any closed form $\eta\in \Lambda^{1,1}(X)$
homologous to $\frac{-2\1 \pi}k c_1(L_0)$, $\nabla$ is the unique 
Chern connection on $L$ with the curvature $\pi^*\eta$.

\hfill

\pstep
Choose a constant metric $h^k$ on $L^{\otimes k}\restrict F=\calo_F$
and let $h$ be its $k$-th root, which is a metric on $L\restrict F$. 
Since  $h^k$ is a constant metric on a trivial line bundle, 
its curvature is flat, and
the Chern connection $\nabla$ associated with $h$ is also flat.

\hfill

{\bf Step 2:} 
The last statement of \ref{_flat_connection_Proposition_}
is clear because for every closed $(1,1)$-form $\eta$ homologous
to $-2\1 \pi c_1(L_0)$ there exists a metric
such that the curvature of its Chern connection
is equal to $\eta$ (\cite{_Grif-Har_}). 
\endproof

\hfill

\definition\label{_local_mono_Definition_}
We define the {\bf fiberwise monodromy},
or {\bf local monodromy},
of $L$ as the monodromy of $\nabla$ in a neighboorhood of a fiber of $\pi$.
Clearly, the fiberwise monodromy is an invariant
of the line bundle $L$. Moreover, $L$ is a pullback
of a line bundle $L'$ on $X$ if and only if the 
fiberwise monodromy of $L$ on each fiber is trivial.

\hfill

\remark
When a fiber $F$ is non-multiple, it is a deformation
retract of its neighbourhood (\cite{_Clemens:degene_,_Persson:degene_}).
When $F$ is a multiple fiber over a general point
of a discriminant, a cyclic covering of its neighourhood 
projects with trivial multiplicity
to a cyclic ramified covering of the base
(\cite[Proposition 3.1]{_Hwang_Oguiso:multiple_}).
In both cases a flat line bundle with trivial local
monodromy is a pullback of a line bundle on the base.

\hfill

\definition\label{_flat_conn_on_restriction_Definition_}
Consider a Lagrangian fibration $\pi:\; M \to \C P^n$,
and let $L$ be a primitive line bundle such that 
$L^{\otimes k}=\pi^* \calo(1)$.
Since $\calo(1)$ is trivial on $\C^n=\C P^n\backslash H$,
\ref{_flat_connection_Proposition_} also gives a flat connection on 
$L\restrict {\pi^{-1}(\C P^n\backslash H)}$.
It can be obtained explicitly as follows.
Let $\nabla$ be the fiberwise flat connection 
defined in \ref{_flat_connection_Proposition_};
its curvature is $\frac 1 k \pi^*\omega_0$, where
$\omega_0$ is the Fubini-Study form on $\C P^n$.
On $\C^n= \C P^n \backslash H$, $\omega_0$ is exact,
$\omega_0= \frac{\1}{2\pi} \6\bar\6 \log\left(1+ \sum_{i=1}^n |z_i|^2\right)$;
let $\theta:= \frac{1}{\pi}d^c \log\left(1+ \sum_{i=1}^n
|z_i|^2\right)$ be its antiderivative. Then
$d(\theta)=\omega_0$. We modify 
$\nabla$ by taking $\nabla_\theta:= \nabla -\frac 1 k \pi^*\theta$,
then $\nabla_\theta^2= \nabla^2- d\theta=0$.

\hfill

\remark\label{_monodromy_then_primitive_Remark_}
If we prove that the monodromy of the flat connection
$\nabla_\theta$ is trivial on $\pi^{-1}(z)$, where $z\in D$
is a general point of the discriminant,
it will follow that it is trivial everywhere. Indeed, 
the fundamental group does not change if
one removes a subvariety of complex codimension 2.
By \ref{_trivial_on_fibers_Theorem_}, the fiberwise monodromy
is trivial on general fibers of $\pi$, 
hence it can theoretically be non-trivial
only on $\pi^{-1}(z)$, where $z\in D$
is a general point of the discriminant.
However, triviality of fiberwise monodromy
of $L$ implies that $L \cong \pi^* \pi_* L$, hence
implies primitivity of $\pi^*\calo(1)$.
Therefore, \ref{_primitivity_Theorem_}
is implied by the following topological result.
Note that we state it for an arbitrary abelian
fibration: the proof is purely topological.

\hfill

\theorem\label{_pi_1_preimage_Theorem_}
Let $\pi:\; M \to \C P^n$ be an abelian 
fibration with no multiple fibers in codimension 1,
and $H \subset \C P^n$ a hyperplane section.
Let $M_0:= \pi^{-1}(\C P^n \backslash H)$.
Then the fundamental group $\pi_1(M_0)$ 
is generated by the smallest normal subgroup containing
$\pi_1(F)$, where $F \subset M_0$ is a general fiber of $\pi$.

\subsection{Fundamental groups of abelian fibrations}

In this subsection we prove 
\ref{_pi_1_preimage_Theorem_}. This also proves
\ref{_primitivity_Theorem_}.

\hfill

{\bf Step 1:} Let $D_0\subset \C^n$ be the discriminant
(that is, the set of critical values) of $\pi:\; M_0 \to \C^n$.
Denote by $X_1$ the complement 
$\C^n \backslash D$ and by $M_1$ its preimage in $M_0$.
By \cite[IX, Cor. 5.6]{_Grothendieck:SGA1_}, 
the natural map $\pi_1(M_1)\to \pi_1(M_0)$ is surjective.
However, the map $\pi:\; M_1 \to X_1$
is a Serre fibration with torus fibers,
giving an exact sequence
\[
\pi_1(F) \to \pi_1 (M_1) \to \pi_1(X_1) \to 0,
\]
where $F$ is a general fiber.
The fundamental group $\pi_1(X_1, p)$  of $X_1$ is generated by
the following loops $\Gamma_z$ starting at $p$. We choose a smooth
point  $z$ in $D$, connect $p$ to $z'$ in a small
neighbourhood of $z$ by a path $\gamma_z$. Then we do a small loop $l_{z}$ 
starting and ending in $z'$ and going around $D$ near $z$ and
return back using $\gamma_z$, obtaining $\Gamma_z=
\gamma_z \circ l_z \circ \gamma_z^{-1}$.

To prove that $\pi_1(M_0)$ is generated by the image
of $\pi_1(F)$, it remains to show that the preimage of each generator
$\Gamma_z$ in $\pi_1(M_1)$ is homotopic in $\pi_1(M_0)$ to
an element of $\pi_1(F)^u$, where $(\cdot )^u$ denotes
the conjugation with $u\in\pi_1(M_0)$.

\hfill

{\bf Step 2:} 
Let $z''$ be a preimage of $z'$ 
in a neighbourhood of a preimage of $z$
which belongs to a reduced component of $\pi^{-1}(D)$,
and $\tilde l_z$ a lift of $l_z$
starting in $z''$, going around
$\pi^{-1}(D)$ and getting back to $z''$.
Clearly, its image in $M_0$ is contractible.

\hfill

{\bf Step 3:}
Since $\pi:\; M_1 \to X_1$
is a Serre fibration with connected fibers,
every path in $X_1$ can be lifted to a path in $M_1$,
and this lift is homotopically unique  up to a 
composition with a loop in the fiber of $\pi$.
Let $\tilde \gamma_z$ be such a lift of $\gamma_z$.
Then, $\Gamma_z$ can be lifted to a path
$\tilde \gamma_z\circ A_1 \circ \tilde l_z \circ A_2  \circ\tilde \gamma_z^{-1}$,
where $A_i$ are paths in $\pi^{-1}(z')\cong F$ and
$\tilde l_z$ the loop around a reduced component of $\pi^{-1}(D)$
defined in Step 2. Since $\tilde l_z$ is
contractible in $\pi_1(M_0)$, we obtain that
$\pi_1(M_0)$ is generated by 
$\tilde \gamma_z\circ A_1 \circ A_2\circ \tilde \gamma_z^{-1}$,
where $A_1, A_2$ are loops in a smooth fiber of $\pi$.
This finishes the proof of 
\ref{_pi_1_preimage_Theorem_}.
\endproof

\subsection{Primitivity and multiple fibers in codimension 1}

\theorem\label{_primitive_multfibers_equiv_Theorem_}
Let $\pi:\; M \arrow \C P^n$ be a Lagrangian fibration
on a hyperk\"ahler manifold of maximal holonomy, 
and $H\subset \C P^n$ a hyperplane section. 
Then the following assertions are equivalent.
\begin{description}
\item[(i)]  The homology class of $\pi^{-1}(H)$ is
primitive.

\item[(ii)] The map $\pi$  has has no multiple fibers
in codimension 1.
\end{description}

\pstep  
The implication (ii) $\Rightarrow$ (i) is proven 
in Section \ref{_primitivity_Section_} using the ETMDPS vanishing theorem.
For the converse implication: let $D_1$ be an irreducible
component of the discriminant $D$; assume it has multiplicity $\mu$.
We need to show that the preimage of the hyperplane section $H$
is not primitive.

\hfill

{\bf Step 2:} Arguing ad absurdum, 
assume that $\pi^{-1}(D_1)$ is multiple of
multiplicity $\mu$, but $[\pi^{-1}(H)]$ is primitive in
$H_{4n-2}(M, \Z)$. 
Clearly, for some $k\in \Z^{>0}$, the divisor $D_1$ is homologous to $k[H]$.
Since $[D_1]= k[H]$, $\pi^{-1}(D_1)$ is the zero set
of a section $s\in H^0(\C P^n, \calo(k))$.
The pullback of this section has multiplicity $\mu$;
since $\pi^*(\calo(1))$ is primitive, $k$ is divisible by $\mu$.
The universal coefficients formula implies that $H^2(M, \Z)$
is torsion-free (\ref{_torsion_free_Claim_}). Therefore, $\pi^*s= s_0^\mu$, where
$s_0\in H^0(M, \pi^*\calo(k/\mu))$. Since $\pi^*\calo(k/\mu)$
is trivial on all fibers of $\pi$, the section $s_0$
is a pullback of a section $s_1 \in H^0(\C P^n,\calo(k/\mu))$.
This is impossible, because $s$, and hence $s_1$, vanishes
on $D_1$, which has degree $k$.
\endproof


\section{Line bundles and multiple fibers}
\label{_line_bundles_Section_}


Let $\pi:\; M \to \C P^n$ be a Lagrangian fibration, and
$L$ a primitive root of the line bundle $\pi^*(\calo(1))$ of degree $d>1$.
The aim of this section is to show that $H^0(L\otimes \pi^*(\calo(i)))=0$
for all $i$.

\subsection{Fundamental group of the preimage of $\C^n$}
\label{_fundamental_preimage_Subsection_}

\definition\label{_essential_multiple_Definition_}
Let $\pi:\; M \to \C P^n$ be a Lagrangian fibration,
$L$ a primitive root of $\pi^*(\calo(1))$, and
$M_0:= \pi^{-1}(\C^n)$ be the preimage of the
complement to a hyperplane. Consider 
the flat connection $\nabla_\theta$ of $L\restrict {M_0}$.
As shown in \ref{_monodromy_then_primitive_Remark_},
$\nabla_\theta$ has trivial monodromy if
and only if $\pi^{-1}(\calo(1))$ is primitive.
The union of all irreducible components of $\pi^{-1}(z)$, 
$z\in \C P^n$, for which the fiberwise connection
is non-trivial, is called {\bf an essential
multiple fiber}, and their union {\bf the
essential multiple divisor.} From \ref{_trivial_on_fibers_Theorem_}
it follows that the essential multiple
divisor is a collection of some irreducible components of 
the preimage of the discriminant.
{\bf The order} of the irreducible component of the essential
multiple divisor is the cardinality
of the image of the monodromy of $\nabla_\theta$
in a neighbourhood of $\pi^{-1}(z)$.

\hfill

\claim\label{_essential_multip_Claim_}
Each component of an 
essential multiple divisor is multiple,
and its order is equal to its multiplicity.

\hfill

\proof
Let $A$ be an irreducible component of the 
essential mutiple divisor of multiplicity $d$.
Then $\calo(dA)= \pi^*(\calo(B))$, where $B=\pi(A)$.
This implies that the bundle $\calo(A)$ admits a section $\nu$
such that its $d$-th 
power $\nu^d$ is a pullback of a section $\nu_1$ of $\calo(B)$.
Since the zeros of $\nu_1$ are of order 1,
the monodromy of $\nu$ along a loop going around $B$ has order $d$.
\endproof

\hfill

Let $L$ be the primitive root of $\pi^*(\calo(1))$, $L^d=\pi^*(\calo(1))$.
In this subsection, we are going to show that $H^0(M,L)=0$.
Otherwise, the zero divisor $D_1$ of a section $\nu\in H^0(M,L)$
is homologous to $\frac 1 d [H]$, where $H$ is the
hyperplane section. Since $L$ is trivial on all
regular fibers of $\pi$, this implies that $\pi(D_1)$
is a union of several irreducible components of the
discriminant. Therefore, $\pi^{-1}(\pi(D_1))$ is a divisor
of multiplicity $d$, and $D_1= \frac 1 d \pi^{-1}(\pi(D_1))$.
Then $\frac 1 d [H]=c_1(L)= \frac 1 d [\pi^{-1}(\pi(D_1))]$.
This implies that
$\pi(D_1)$ is homologous to a hyperplane, 
hence $\pi(D_1)$ {\em is} a hyperplane $H\subset \C P^n$,
its preimage has multiplicity $d$, and
$L=\calo(\frac 1 d \pi^{-1}(H))$. We are going to show
that this is impossible, {\em ipso facto} $H^0(M,L)=0$.

\hfill

The idea of the proof is
simple and geometric.
Consider an abelian fibration
$\pi:\; M \to X$, which is multiple over a divisor $Z\subset X$.
Locally in $X$ in a neighbourhood of a general point
of the discriminant, there is always a finite ramified
covering  $X'\to X$ such that the pullback of
$\pi$ to $X'$ is a fibration without multiple
fibers in codimension 1 (\cite{_CKV:min_mult_}).
In the situation we discussed, we are able to produce
such a ramified covering globally on $M$.
Also we have $X=\C P^n$ and
the ramification divisor is a hyperplane.
Since $\pi_1(\C P^n \backslash H)=0$,
this is impossible, and hence the essential
multiple divisor of a Lagrangian fibration
$\pi:\; M \to \C P^n$ cannot be a preimage
of a hyperplane. This was a sketch of a proof of
\ref{_complement_to_C^n_Theorem_} below.

\hfill

\theorem\label{_complement_to_C^n_Theorem_}
Let $\pi:\; M \to \C P^n$ be a Lagrangian
fibration, and $L$ a primitive root $\pi^*(\calo(1))$,
$L^d=\pi^*(\calo(1))$. Then $H^0(M,L)=0$.

\hfill

\pstep
As we explained earlier in this subsection,
$H^0(M,L)\neq 0$ only if a hyperplane 
$H\subset \C P^n$ is an irreducible component
of the discriminant of $D$ and $\pi^{-1}(H)$
is a divisor of multiplicity $d$. Denote the
 divisor $\frac 1 d \pi^{-1}(H)$ by $H'$.
Clearly, $\dim H^0(M,L)=1$, because
all sections of $L$ vanish on $\pi^{-1}(H)$.
Consider a section of $L$ as a subvariety
$M'$ of a total space $\Tot(L)$.
Clearly, $M'$ is a ramified cover of $M$
of degree $d$, and $H'$ is its ramification divisor.
Consider the Stein factorization of the
projection $M' \to \C P^n$, with
$\pi':\; M' \to B'$ being a fibration
with connected fibers, and $\sigma:\; B'\to \C P^n$
finite. We are going to show that
$\pi':\; M' \to B'$ is an abelian
fibration with no multiple fibers
over a general point of $H$.

\hfill

{\bf Step 2:} Let $x\in H'\subset M$ be a general point.
Locally around $H'$, the map $\pi$ is written in coordinates as
$(p_1, ..., p_n, q_1, ... q_n)\mapsto (p_1^d, p_2, ..., p_n)$
(\cite{_Hwang_Oguiso:multiple_,_CKV:min_mult_}), where
$p_1=0$ is the equation defining $H'$.
Since $\pi^*(\calo(1))$ has zero of order $d$ in $H'$,
the section of $L$ has zero of order 1 in $H'$.
Therefore, after the local trivialization, the global
section of $L$ has the same zeros as the function $p_1$.
This implies that the natural map $M'\to B'$
takes $(p_1, ..., p_n, q_1, ... q_n)$ to $(p_1, ..., p_n)$,
and it is unramified in a general point of $H'$:
we proved the claim made in the end of Step 1.

\hfill

{\bf Step 3:} Since $\pi_1(\C P^n\backslash H)=0$,
the ramified covering $B'\to \C P^n$ has trivial ramification index.
However, for a general point $x\in H$, the fiberwise
monodromy group of $L$ acts freely on its preimage
in $M'$, because $M'\to \C P^n$ has no multiple fibers
in codimension 1. Therefore, the map $M'\to M$
is a covering outside of codimension 2; this is
impossible because $\pi_1(M)=0$.
This finishes the proof of 
\ref{_complement_to_C^n_Theorem_}.
\endproof

\subsection{Irreducible components of the 
essential multiple divisor}
\label{_essential_irreducible_Subsection_}

The aim of this subsection is the following

\hfill

\proposition\label{_multiplicities_coprime_Proposition_}
Let $\pi:\; M \to \C P^n$ be a Lagrangian fibration,
and $D_1, ..., D_k\subset M$ are the irreducible
components of the essential multiple divisor
(\ref{_essential_multiple_Definition_}).
Let $\mu_i$ be the multiplicity of $D_i$.
Then the numbers $\mu_i$ are pairwise coprime.

\hfill

\pstep
Consider the group $\Lambda:= \langle L\rangle \subset \Pic(M)$ generated by 
$L$, which is defined as the primitive root of $\pi^*\calo(1)$.
Then $\calo(D_i)\in \Lambda$, because $\mu_i D_i$ is the
preimage of a divisor in $\C P^n$ (\ref{_essential_multip_Claim_}).
Let $\Lambda_0:= \frac{\langle L\rangle}{\pi^*\calo(1)}$;
by \ref{_torsion_free_Claim_}, this group is cyclic. 
Consider its subgroup $\Lambda_1$ generated by $\calo(D_i)$, 
$i=1, ..., k$. Clearly, this subgroup is also cyclic.
Recall that the Chinese Remainder Theorem implies
that $\Z/{a\Z} \oplus \Z/{b \Z}$ can be realized
as a subgroup in a cyclic group if and only if $a, b$ are coprime.
Therefore, to prove that $\mu_i$ are pairwise coprime, it
suffices to show that 
$\Lambda_0$ is generated by $\calo(D_1)$, ... $\calo(D_k)$,
and relations $\calo(D_i)^{\mu_i}=0$.

\hfill

{\bf Step 2:}
Let $\Lambda_2\subset \Pic(M)$ be the subgroup
generated by $\calo(D_1)$, ... $\calo(D_k)$.
All the relations in $\Lambda_2$ are obtained from 
$Z_1 \sim Z_2$ whenever $Z_1$ and $Z_2$
are zero and pole divisors of a rational
function $f$.  Since the line bundles $\calo_{D_i}$ are trivial on
fibers of $\pi$, these functions
are constant on the fibers of $\pi$,
hence we can assume that $Z_1$ and $Z_2$
are supported on the preimages of divisors on $\C P^n$.
However, after passing from $\Lambda_2$
to $\Lambda_0:= \frac{\Lambda_2}{\langle\pi^*\calo(1)\rangle}$,
the divisors $Z_i$ vanish, hence the only
relations in $\Lambda_2$ are obtained from
taking a power of $\calo(D_i)$ and
obtaining $\pi^*(\calo(Z))= \calo(D_i)^{t_i}$;
by definition of multiplicity, this gives $t_i\mid \mu_i$.
\endproof

\subsection{Complex line bundle and torsion in cohomology}
\label{_torsion_line_bundle_Subsection_}
 
Here and in the sequel,
``flat bundle'' means a bundle equipped with a flat 
unitary connection.
The following theorem seems to be known, but
we did not find a reference.

\hfill

\theorem\label{_monodromy_on_torsion_Theorem_}
Let $E$ be a flat complex line bundle
on a manifold $X$. Then $E$ is topologically trivial
if and only if its monodromy, restricted
to the torsion part in $H_1(X,\Z)$, is trivial.

\hfill

\pstep
Let $\Pic_{C^\infty}(X)$ be the group
of smooth complex line bundles on $X$.
The exponential exact sequence immediately
implies that an element in $\Pic_{C^\infty}(X)$
is trivial if and only if its first Chern class
is trivial; since $E$ is flat, $c_1(E)$ is torsion.
Let $C(X)$ be the group of characters on
$\pi_1(X)$, that is, the group of homomorphisms
$\pi_1(X) \to U(1)$. We are interested
in the kernel of the natural forgetful map 
$\Phi:\; C(X)\to \Pic_{C^\infty}(X)$.

\hfill

{\bf Step 2:} The kernel of $\Phi$ consists
of unitary complex local systems with non-trivial monodromy
on a trivial vector bundle, that is, the
maps $\pi_1(X) \to \pi_1(S^1)=\Z$.
Clearly, two local systems
$a, b:\; \pi_1(X) \to U(1)$
are topologically isomorphic  whenever
$a b^{-1}$ corresponds to a trivial
line bundle, that is, $a b^{-1}$ is trivial
on the torsion subgroup in $H_1(X, \Z)$. 
The restriction of a character
$u\in C(X)$ to the torsion-free part does not
matter for its image in $\Pic_{C^\infty}(X)$, because any homomorphism
$\Z^n \to U(1)$ can be lifted
to a homomorphism $Z^n \to \Z$,
hence corresponds to a topologically
trivial line bundle. \endproof

\subsection{Sections of line bundles and monodromy}
\label{_sections_monodromy_Subsection_}

The aim of this section is the following

\hfill

\theorem\label{_sections_of_L_otimes_pullback_Theorem_}
Let $\pi:\; M \to \C P^n$ be a Lagrangian
fibration, and $L$ a primitive root of $\pi^*\calo(1)$,
of degree $d>1$. Then $H^0(M, L\otimes \pi^*(L_1))=0$
for all line bundles $L_1$ on $\C P^n$.

\hfill

\pstep
Let $Z\subset \C P^n$ be an irreducible divisor in 
the complex projective space. It is well known
that $H_1(\C P^n\backslash Z)$ is cyclic
(\cite[Chapter 4, Proposition 1.3]{_Dimca:hypersurfaces_}).
By \ref{_monodromy_on_torsion_Theorem_}, this implies
that no trivial line bundle on $\C P^n\backslash Z$ 
can have a flat connection with non-trivial monodromy.

\hfill

{\bf Step 2:}
Suppose that $H^0(M, L\otimes \pi^*(L_1))\neq 0$.
Since $L\otimes \pi^*(L_1)$ is trivial on general fibers,
the zero divisors of its sections are 
of the form $\pi^{-1}(Z)+ \sum a_i D_i$, where
$Z$ is a divisor on $\C P^n$, and $D_i$
are the irreducible components
of the essentially multiple divisor.
Tensoring with $\pi^*(\calo(-Z))$, we obtain
a section of $L\otimes \pi^*(L_1)\otimes \pi^*(\calo(-Z))$
which does not vanish outside of the essentially
multiple divisor. This implies that the line bundle
$L\otimes \pi^*(L_1)\otimes \pi^*(\calo(-Z))$
is trivial on the complement to the essentially multiple divisor.

\hfill

{\bf Step 3:}
Let $B:= \pi(D)$, where $D$ is the essentially multiple divisor.
Here we prove that $B$  cannot
be irreducible if $H^0(M, L\otimes \pi^*(L_1))\neq 0$.
Clearly, $D= \frac 1 \mu \pi^{-1}(B)$, 
where $\mu$ is its multiplicity.
 Clearly,  $L\otimes \pi^*(\calo(-Z))$
is trivial on the complement to the essentially multiple
divisor, because $\calo(i)$ is trivial on 
 the complement to any divisor in $\C P^n$.
Since $L$ has flat connection on  $\C P^n\backslash B$,
it is trivial on $\C P^n\backslash B$ (Step 1). 
This implies that $L=\calo(l D)$, for some integer $l$,
hence $\deg B=1$ and $L= \calo(\frac 1 \mu \pi^{-1}(B))$.
The latter is impossible by
\ref{_complement_to_C^n_Theorem_}.

\hfill

{\bf Step 4:}
Using Step 3, we obtain that it suffices to prove
\ref{_sections_of_L_otimes_pullback_Theorem_}
when the essentially multiple divisor is not irreducible.
Let $D_1, ..., D_k\subset M$ be the irreducible components
of the essentially multiple divisor, $\mu_i$ the
corresponding multiplicities, and $B_i:=\pi(D_i)$.
Since $\mu_i$ are coprime 
(\ref{_multiplicities_coprime_Proposition_}), for $d:= \mu_2 \mu_3...\mu_k$ 
the bundle $L^d$ is trivial on all $D_i$ except $D_1$;
moreover, this bundle has non-trivial local
monodromy over $\pi^{-1}(z)$, for a general $z\in D_1$.
By the same argument which proves 
Step 2, $H^0(M, L\otimes \pi^*(L_1))\neq 0$
implies that for some line bundle $L_2$ on $\C P^n$
the line bundle $L^d\otimes \pi^*(L_2)$ is trivial
on the general fibers over $B_2, ..., B_k$ and
on the fibers over the complement of $\bigcup B_i$.
Therefore, the bundle $L^d\otimes \pi^*(L_2)$
is trivialized by the section we constructed 
on the complement to $D_1$. 
Since $\pi^*(\calo(1))\restrict{\C  P^n\backslash D_1}$
is trivial, this implies that 
$L^d\restrict{\C  P^n\backslash D_1}$ is trivial,
hence $L\restrict{\C  P^n\backslash D_1}$
is a flat line bundle.

\hfill

{\bf Step 5:} 
Since $d$ is coprime with $\mu_1$, triviality of
$L^d \restrict{\C  P^n\backslash D_1}$ implies 
triviality of $L\restrict{\C  P^n\backslash D_1}$.
Indeed, if $dd'=1 \mod \mu_1$, we have
$L^{d'}\restrict{\C  P^n\backslash D_1}\cong L\restrict{\C  P^n\backslash D_1}$.
However, 
since $D_1$ is irreducible, any line bundle which
is trivialized on its complement is isomorphic to
$\calo(l D_1)$, for some integer $l$.
Then $L=\calo(l D_1)$, which is
impossible because $L$ has non-trivial
monodromy over the general points of $D_2, D_3, ..., D_k$.
\endproof


\section{Differential forms with coefficients in line bundles and stability}


%
%
%
%
%

\subsection{Stability of $\Omega^iM$ and its implications}

\definition
Let $F$ be a torsion-free coherent sheaf on $M$.
Define {\bf the degree}
$\deg_\omega F:=\int_M c_1(F) \wedge \omega^{n-1}$,
where $\omega$ is a K\"ahler form. Let
$\slope(F):=\frac{\deg_\omega F}{\text{rank}(F)}$.
A torsion-free sheaf $F$ is called {\bf stable}
if for all subsheaves $F'\subset F$ one has
$\text{slope}(F')<\text{slope}(F)$. If $F$ is  
a direct sum of stable sheaves of the same slope, 
$F$ is called {\bf polystable}.

\hfill

\theorem
On a hyperk\"ahler manifold of maximal holonomy, 
the bundle $\Omega^{2i+1}(M)$ is stable for all
$i$, and the bundle $\Omega^{2i}(M)$
is polystable.  Moreover, $\Omega^{2i}(M)$
is a direct sum of a stable bundle and $\C \Omega^{i}$,
where $\Omega$ denotes the holomorphically symplectic form.

\hfill

\proof
Follows from the Kobayashi-Hitchin correspondence
(\cite{_Donaldson:surfa_,_UY_,_Lubke_Teleman:Book_}).
By Calabi-Yau theorem, $TM$ and its tensor powers
admit a Yang-Mills metric, and a Yang-Mills bundle
is polystable, and stable if the holonomy of the
Yang-Mills connection is irreducible. However,
the holonomy of the Yang-Mills connection on
a hyperk\"ahler manifold of maximal holonomy
is $Sp(n)$, and the irreducible decomposition
of the Grassmann algebra for the fundamental
representation of $Sp(n)$ is described in 
\cite{_Weyl:invariants_} (see also 
\cite{_Howe:duality_,_Verbitsky:balanced_,_Verbitsky:SYZ_}).
\endproof

\hfill

\remark Clearly, a tensor product of a stable bundle and
a line bundle is also stable.

\hfill

The following corollary is interesting,
but we are not going to use it further in this paper.

\hfill

\corollary\label{_L_valued_form_non_zero_on_D_Corollary_}
Let $\pi:\; M \rightarrow \C P^n$ be a
Lagrangian fibration, and
$L$ the primitive root of $\pi^*(\calo(1))$, such that
$L^d=\pi^*(\calo(1))$. Let $u$ a section of 
$L\otimes \Omega^*(M)$.
Then $u$ is non-zero outside of a codimension 2 subvariety.

\hfill

\proof
Since $L$ is trivial on any general fiber $F$ of $\pi$
(\ref{_trivial_on_fibers_Theorem_}) the  zero divisor of $u$ is 
the preimage of a divisor $Z\subset \C P^n$.
Indeed, $TM\restrict F$ is an extension of trivial bundles,
and any section of its tensor powers is zero or vanishes everywhere.

This implies that the codimension 1 part of the zero set of $u$
is  $\pi^{-1}(Z)$. Then the degree
of the line sub-bundle $V$ of $L\otimes \Omega^{*}(M)$
generated by $u$ is at least 
$\deg \pi^{-1}(Z)= \deg_{\C P^n} Z\cdot  \deg_\omega \pi^*(\calo(1))$.
We rescale the K\"ahler form on $M$ in such a way
that $\int_{\pi^{-1}(H)}\omega^{2n-1}=1$, where $H$ is a
hyperplane divisor in $\C P^n$.
This is equivalent to $\int_{\pi^{-1}(Z)}\omega^{2n-1}=\deg Z$
for any divisor $Z\subset \C P^n$. Since $L^{\otimes k}=\pi^*\calo(1)$
and $\deg_\omega \pi^*\calo(1)= \int_{\pi^{-1}(H)}\omega^{2n-1}=1$, we have
$\deg_\omega L= 1/d$.
Then 
\[ 
\slope(V)\geq \int_{\pi^{-1}(Z)}\omega^{2n-1} =
\deg Z > \slope(L\otimes \Omega^{2i}(M))=
\deg_\omega L = d^{-1}.
\]
This contradicts the polystability of $L\otimes \Omega^{*}(M)$.
\endproof

\hfill

\remark
From \ref{_n+1+H_0_Omega^iotimesL_Corollary_}, 
we obtain that $\dim H^0(\Omega^* M \otimes L)\geq n+1$,
hence the non-zero holomorphic $L$-valued forms exist.
\ref{_L_valued_form_non_zero_on_D_Corollary_} 
implies that these forms might vanish
only in codimension $\geq 2$ sets.

\subsection{Differential forms with coefficients in line bundles}

In the rest of this section, we prove  \ref{_main_intro_Theorem_}
which claims that a Lagrangian fibration $\pi:\; M \to \C P^n$ 
has no multiple fibers in codimension 1. We prove it
by showing that the bundle $\pi^*(\calo(1))$ is primitive.

\hfill

We deduce the primitivity of $\pi^*(\calo(1))$ from 
\ref{_n+1+H_0_Omega^iotimesL_Corollary_}, which
claims that $H^0(L\otimes \Omega^*(M))\neq 0$
whenever  $L^k=\pi^*(\calo(1))$.
Using an $L$-valued holomorphic differential
form, we construct a section of $\pi^* L_1 \otimes L$,
coming into a contradiction with 
\ref{_sections_of_L_otimes_pullback_Theorem_}.

\hfill

\proposition\label{_diff_forms_to_line_bundle_Proposition_}
Let $\pi:\; M \to \C P^n$ be a Lagrangian fibration
on a hyperk\"ahler manifold, $L^k=\pi^*(\calo(1))$
a primitive line bundle, $k >1$. Assume that 
$H^0(L\otimes \Omega^*(M))\neq 0$.
Then there exists a non-zero section on $\pi^*L_1\otimes L$,
where $L_1:= \pi^*(K_{\C P^n}^{-1})$.

\hfill

\pstep
Let ${\cal F}_1 \subset \Omega^*(M)$ be the sheaf of all
differential forms $\eta$ which are divisible by $\pi^* \Lambda^n \C P^n$,
that is, satisfy $\eta\wedge \pi^*\theta=0$
for any 1-form $\theta\in \Lambda^{1,0}(\C P^n)$.
Starting with an $L$-valued differential form $\nu\in H^0(L\otimes \Omega^*(M))$,
we are going to produce a non-zero 
form $\nu_1\in H^0(L\otimes {\cal F}_1\otimes \pi^*(K_{\C P^n})^{-1})$.

Clearly, the sheaf $\Omega^* \C P^n \otimes K_{\C P^n}^{-1}= \Lambda^* T\C P^n$
is globally generated. We are going to multiply $\nu$ by a section
of 
$\pi^*(\Omega^* \C P^n \otimes K_{\C P^n}^{-1})\subset \Omega^* M \otimes\pi^*(K_{\C P^n}^{-1})$
to obtain a non-zero section of ${\cal F}_1\otimes \pi^*(K_{\C P^n}^{-1})$.

Let $dp_i, dq_i$ be a local
frame in a regular point $x\in M$ of $\pi$ such that
$p_i$ are coordinates on $\C P^n$ and the holomorphically
symplectic form is $\Omega=\sum_i dp_i \wedge dq_i$.
Then ${\cal F}_1$ in a neighbourhood of $x$ are differential
forms divisible by $dp_1\wedge ... \wedge dp_n$.
If $\nu$ is locally written as a sum of monomials
on $dp_i$, $dq_i$, which contains a monomial
$dp_{i_1}\wedge ... \wedge dp_{i_a}\wedge dq_{j_1}\wedge ... \wedge dq_{j_b}$,
we can always multiply it by a complimentary monomial on $dp_i$
to obtain $dp_1\wedge ... \wedge dp_n \wedge dq_{j_1}\wedge ... \wedge dq_{j_b}$.
This construction can be applied globally, because 
$\pi^*(\Omega^* \C P^n \otimes K_{\C P^n}^*)$ is globally generated.
Therefore, in the assumptions of 
\ref{_diff_forms_to_line_bundle_Proposition_}
we can additionally assume that $\nu$ is a section
of $L \otimes {\cal F}_1 \otimes \pi^* L_1\subset L\otimes \Omega^*(M)$.

\hfill

{\bf Step 2:} Let $\nu$ be a non-zero 
section of $L\otimes {\cal F}_1 \otimes \pi^* L_1$.
Denote by $T_\pi M\subset TM$ the sheaf of $\pi$-vertical tangent vectors,
that is, the kernel of $d\pi:\; T_y M\to T_{\pi(y)} \C P^n$,
and let $\Lambda^* T_\pi M\subset \Lambda^* T M$ be the
corresponding sheaf of polyvectors.
We are going to contract $\nu$ with a global holomorphic section of 
$\Lambda^* T M \otimes \pi^* K^{-1}_{\C P^n}$
to obtain a non-zero section of 
$L\otimes \pi^*\Omega^n \C P^n\otimes \pi^*K^{-1}_{\C P^n}\otimes \pi^* L_1=L\otimes \pi^* L_1$,
thus finishing the proof.

Using the same coordinate system as in Step 1, we can
locally express
the holomorphically symplectic form $\Omega$
as $\sum_i dp_i \wedge dq_i$, and its dual
bivector as $\sum_i\frac d{dp_i} \wedge \frac d{dq_i}$.
This bivector defines
an isomorphism $\Omega^1 M \to TM$,
taking $dp_i$ to $\frac d{dq_i}$.
Extending this isomorphism to 
$\pi^*\Omega^i \C P^n=\Lambda^*(\langle dp_i\rangle)$
we obtain a map $\pi^*\Omega^i \C P^n\to \Lambda^i T_\pi M$,
defining an isomorphism in all regular points of $\pi$.
Since the sheaf $\Omega^i \C P^n\otimes  K^{-1}_{\C P^n}$
is globally generated, the sheaf $\Lambda^i T_\pi M\otimes  \pi^*K^{-1}_{\C P^n}$
is globally generated on the set of regular points of $\pi$.

Consider again a monomial decomposition of
$\nu \in H^0(L\otimes {\cal F}_1 \otimes \pi^* L_1)$, and assume that
it contains a monomial $dp_1\wedge ... dp_n \wedge dq_{j_1}\wedge ... dq_{j_b}$
(with non-zero coefficient). Contracting $\nu$
with the polyvector $\frac d{q_{j_1}}\wedge ... \frac d {q_{j_b}}$,
we obtain a non-zero section of $L\otimes \pi^*\Omega^n \C P^n\otimes \pi^* L_1$.
On the other hand, the sheaf $\Lambda^i T_\pi M\otimes  \pi^*K^{-1}_{\C P^n}$
is globally generated, hence the contraction
\[
L\otimes {\cal F}_1 \otimes 
\pi^* L_1\times \Lambda^i T_\pi M\otimes  \pi^*K^{-1}_{\C P^n}
\to L\otimes \pi^*\Omega^n \C P^n\otimes  \pi^*K^{-1}_{\C P^n}\otimes\pi^* L_1
\]
takes $\nu$ to a non-zero global section of 
$L\otimes \pi^*\Omega^n \C P^n\otimes  \pi^*K^{-1}_{\C P^n}\otimes\pi^* L_1=L \otimes\pi^* L_1$.
\endproof

\hfill



We can finish the proof of our main theorem now.

\hfill

\theorem \label{_no_mult_main_end_Theorem_}
Let $\pi:\; M \to \C P^n$ be a
Lagrangian fibration on a hyperk\"ahler manifold.
Then $\pi$ has no multiple fibers
in codimension 1.

\hfill

\proof
By \ref{_primitive_multfibers_equiv_Theorem_}, 
it suffices to show that the line bundle
$\pi^*\calo(1)$ is primitive. Arguing
ad absurdum, assume that $L^{\otimes \mu} = \pi^*\calo(1)$,
where $\mu>1$. By \ref{_n+1+H_0_Omega^iotimesL_Corollary_}, 
the bundle $\Omega^* M \otimes L$
has non-zero sections. 
By \ref{_diff_forms_to_line_bundle_Proposition_},
this implies that $L\otimes \pi^* L_1$
has non-zero sections for some line bundle $L_1$ on $\C P^n$.
However, 
this is impossible by \ref{_sections_of_L_otimes_pullback_Theorem_}.
\endproof

\hfill

{\bf Acknowledgments.} We are very grateful to Frederic Campana,
Jason Starr, Andrey Soldatenkov and   Yoon-Joo Kim for interesting
discussions and advice. When we were finishing this paper,
Evgeny Shinder contacted us to let us know about the
paper \cite{_KOS:obstructions_} independently proving similar results.
We are grateful to the authors of \cite{_KOS:obstructions_} for 
letting us know about their work. 

\hfill

\noindent {\sc Ljudmila Kamenova\\
Department of Mathematics, 3-115 \\
Stony Brook University \\
Stony Brook, NY 11794-3651, USA,} \\
\tt kamenova@math.sunysb.edu
\\

\noindent {\sc Misha Verbitsky\\
\tt verbit2000@gmail.com
}

\end{document}